\documentclass{article}
\usepackage{indentfirst,latexsym,bm,amsmath,amsthm,amsfonts,array,amssymb}
\begin{document}
\title{Toric Rigid Spaces}
\author{BinYong Hsie
      \\{\small{LMAM, Department of Mathematics, PeKing University, BeiJing,
      P.R.China,
      100871 }}
      \\{\small E-mail: byhsie@math.pku.edu.cn}
      \\ZhiBin Liang
      \\{\small{LMAM, Department of Mathematics, PeKing University, BeiJing,
      P.R.China,
      100871 }}
      \\{\small E-mail: lzb\_1979@hotmail.com}}
\date{}
\maketitle{}

\begin{abstract}
This paper gives a method to construct rigid spaces, which is
similar to the method used to construct toric schemes.
\end{abstract}

{\small{{\bf Keywords:} cone, fan, strongly convex rational
polyhedral cone, toric affinoid algebra, toric rigid space, toric
scheme

MR(2000) Subject Classification: 11E95, 11G25

\section{Introduction}
Toric geometry provides an important way to see many examples and
phenomena in algebraic geometry. Toric geometry studies toric
schemes, which are a special type of schemes. At least, they are
all rational. However they provide intuition for us to study
algebraic geometry and they make computation much easier.

Toric schemes correspond to objects like the simplicial complexes
studied in algebraic topology. To construct a toric scheme, we
need a fan, which consists of (strongly convex rational
polyhedral) cones. If a cone is contained in this fan, so are its
faces. A strongly convex rational polyhedral cone is, by
definition, a cone with an apex at the origin, generated by a
finite number of vectors, and satisfies some conditions.  Such a
cone $\sigma$ provides a finite generates semigroup and this
semigroup gives a neotherian (semigroup) algebra
$\mathbb{Q}[S_{\sigma}]$ over $\mathbb{Q}$. Therefore each cone
provides an affine scheme. Patching these schemes in some way, we
get a separated integral scheme over $\mathbb{Q}$. Call it the
toric scheme associated to the given fan.

Luckily, the authors find this method can be used to construct
rigid spaces, too. Taking the $p$-adic completion of
$\mathbb{Q}_{p}[S_{\sigma}]$, we get an affinoid algebra and then
an affinoid rigid space. Patching all these affinoid rigid spaces,
we get a rigid space, and call it the toric rigid space associated
to the given fan. We find that its natural reduction is exactly
the toric scheme over $\mathbb{F}_{p}$ associated to the same fan.

There are many problems occurring. For example, what is the
relation between the toric rigid space that we get and the rigid
analytification of the toric scheme defined over $\mathbb{Q}_{p}$
associated to the same fan? In general, they are not the same, but
sometimes they are. How can we understand such a phenomenon? We
shall point out that the toric rigid spaces that we construct are
more natural than the latter.

Professor S. Bosch and B. Le Stum gave several courses on rigid
geometry at Peking university. Thanks to their work, we develop
our ideas in this paper.
\section{Toric affinoid algebras and toric rigid spaces}
\subsection{Construction of toric rigid spaces}
Let $K$ be a field with a complete non Archimedean absolute value
and $\bar{K}$ be its algebraic closure. Assume its residue field
$k$ is of characteristic $p$.

 Let $\mathrm{N}$ be a lattice which is isomorphic to
$\mathbb{Z}^n$ for some positive integer $n$ and $\sigma$ be a
strongly convex rational polyhedral cone in the vector space
$\mathrm{N}_{\mathbb{R}}=\mathrm{N}\bigotimes_{\mathbb{Z}}\mathbb{R}.$
 A strongly convex rational polyhedral cone is a cone with an apex at the
origin, generated by a finite number of vectors; ``rational" means
that it is generated by vectors in the lattice $\mathrm{N}$, and
``strong" means that it contains no line through the origin.

Write $\mathrm{V}$ for $\mathrm{N}_{\mathbb{R}}$, then the dual
cone of $\sigma$ is defined as
\begin{equation}
\sigma^{\vee}:\;=\{u\in
\mathrm{V}^{*}:\;<u,v>\geqslant0,\:\forall\; v\in\sigma\}.
\end{equation}
We denote $\mathrm{M}=Hom(\mathrm{N},\mathbb{Z})$ as its dual
lattice with dual pairing denoted by $<\:\;,\:\; >$. Set
\begin{equation}
S_{\sigma}:\;=\sigma^{\vee}\cap\mathrm{M},
\end{equation}
which is a semigroup. As a result, $S_{\sigma}$ is finitely
generated.

 Given such an $S_{\sigma}$, we define an  affinoid algebra
 $K\langle S_{\sigma}\rangle$ associated to it.
 $K\langle S_{\sigma}\rangle$ is just the $p$-adic completion of the (semigroup) algebra $K[S_{\sigma}]$.
 It's equivalent to say that
 \begin{equation}
K\langle S_{\sigma}\rangle:=K\bigotimes _{R}\lim_{\leftarrow}
R/p^{n}R [S_{\sigma}]. \end{equation} Here, $R$ is the integral
ring of $K$. It is easy to check that it is a finitely generated
commutative affinoid algebra. We call it the toric affinoid
algebra associated to $\sigma .$

We denote $\{\chi^{u}|u\in S_{\sigma}\} $  as a basis of the
vector space $ K\langle S_{\sigma}\rangle$, and we have
$$\chi^{u_{1}}\cdot\chi^{u_{2}}=\chi^{u_{1}+u_{2}}.$$
The unit $1$ is just $\chi^{0}$. Generators $\{u_{i}\}$ for the
semigroup $S$ determine generators $\chi^{u_{i}}$ for the
$K$-algebra $K\langle S\rangle$. Denote $X_{i}$ for $\chi^{u_{i}}$
for simple, then
$$K\langle S\rangle=K\langle X_{1},X_{2},\cdot\cdot\cdot,X_{n}\rangle=K\langle\zeta_{1},\zeta_{2},\cdot\cdot\cdot,\zeta_{n}\rangle/\mathrm{I},$$
where $n$ is the number of the generators of $K$-algebra $K\langle
S\rangle$, each $\zeta_{i}$ is an indeterminate element, and
$\mathrm{I} $ is an ideal of
$K\langle\zeta_{1},\zeta_{2},\cdot\cdot\cdot,\zeta_{n}\rangle$.

Set
\begin{equation}
\mathrm{U}_{\sigma}=\mathrm{Sp}(K\langle S_{\sigma}\rangle)
\end{equation}
as the corresponding  rigid space.

Now we begin to consider a fan $\Delta$ in
$\mathrm{N}_{\mathbb{R}}$ and the rigid space defined according to
it. A fan $\Delta$ is defined as the collection of ``strongly
convex rational polyhedral cones" satisfying the following two
conditions:
\begin{quote}
1. Every face of a cone in $\Delta$ is a cone in $\Delta$.

2. The intersection of two cones in $\Delta$ is a face of each
other.
\end{quote}

When $\tau$ is a face of $\sigma$ in $\Delta$, $S_{\sigma}$ is
contained in $S_{\tau}$.We can find a $u\in S_{\sigma}$ such that
$\tau=\sigma\bigcap u^{\perp}$ and
$$S_{\tau}=S_{\sigma}+\mathbb{Z}_{\geq 0}(-u).$$
Then we have
$$K\langle S_{\tau}\rangle=(K\langle S_{\sigma}\rangle)_{\chi^{u}}.$$
 So we get an open embedding
$i_{\sigma}^{\tau}:\;\mathrm{U}_{\tau}\hookrightarrow
\mathrm{U}_{\sigma}$. We glue all these $\mathrm{U}_{\sigma}$ (
$\sigma\in \Delta$) by the following way. If $\sigma,
\sigma_{1}\in \Delta$ and $ \sigma \cap \sigma_{1} =\tau$, we glue
$\mathrm{U}_{\sigma}$ and $\mathrm{U}_{\sigma_{1}}$ by the
immersions $i_{\sigma}^{\tau}$ and $i_{\sigma_{1}}^{\tau}$. We
denote the rigid space we get by $\mathrm{R}(\Delta)$. Write
$\mathrm{T}(\Delta)$ for the corresponding toric scheme in the
same sense as in \cite{Fulton}.

Now we come to the following proposition.

{\bf{Proposition 1.}} $\mathrm{R}(\Delta)$ is a rigid space.


We want to show that all toric rigid spaces are separated. We need
the following lemma.

{\bf Lemma 1.} If $\sigma$ and $\tau$ are two strongly convex
polyhedral cones in  $\mathrm{N}_{\mathbb{R}}$, we have
$\sigma^{\vee}+\tau^{\vee}=(\sigma\cap\tau)^{\vee}$.
\begin{proof}
First we prove $\sigma^{\vee}+\tau^{\vee}\subseteq
(\sigma\cap\tau)^{\vee}$. If $u\in\sigma^{\vee}$ and
$v\in\tau^{\vee}$, and $t\in\sigma\cap\tau $, then
$<u,t>+<v,t>\geq0$ by the definition. So $u+v\in
(\sigma\cap\tau)^{\vee}.$

Then we prove
$\sigma^{\vee}+\tau^{\vee}\supseteq(\sigma\cap\tau)^{\vee}$. It is
equivalent to show that
$(\sigma^{\vee}+\tau^{\vee})^{\vee}\subseteq((\sigma\cap\tau)^{\vee})^{\vee}=\sigma\cap\tau$.

From $\sigma^{\vee}\subseteq\sigma^{\vee}+\tau^{\vee}$, we get
$\sigma=(\sigma^{\vee})^{\vee}\supseteq(\sigma^{\vee}+\tau^{\vee})^{\vee}.$
And from $\tau^{\vee}\subseteq\sigma^{\vee}+\tau^{\vee}$, we get
$\tau=(\tau^{\vee})^{\vee}\supseteq(\sigma^{\vee}+\tau^{\vee})^{\vee}.$
Therefore,
$\sigma\cap\tau\supseteq(\sigma^{\vee}+\tau^{\vee})^{\vee}$.
\end{proof}

From the proof of the above lemma, we can find  it is also right
that
\begin{equation}
(\sigma^{\vee}+\tau^{\vee})\bigcap\mathrm{M}=((\sigma\cap\tau)^{\vee})\bigcap\mathrm{M}.\label{sep}
\end{equation}
Now we have sufficient preparations to prove the following
theorem.

{\bf Theorem 1.} $\mathrm{R}(\Delta)$ is separated for a given fan
$\Delta$ in $\mathrm{N}_{\mathbb{R}}$.
\begin{proof}
$\{\mathrm{U}_{\sigma}:\;\sigma$ is a cone in $ \Delta\}$ is an
affinoid covering of $\mathrm{R}(\Delta)$. We have to prove two
facts that $\mathrm{U}_{\sigma}\bigcap\mathrm{U}_{\tau}$ is
affinoid for $\sigma,\tau$ in $\Delta$, and that
$\mathcal{O}(\mathrm{U}_{\sigma}\bigcap\mathrm{U}_{\tau})$ is
generated by the canonical image of
$\mathcal{O}(\mathrm{U}_{\sigma})$ and
$\mathcal{O}(\mathrm{U}_{\sigma})$.

Because
$\mathrm{U}_{\sigma}\bigcap\mathrm{U}_{\tau}=\mathrm{U}_{\sigma\bigcap\tau}$,
we have got the first fact.

As $\mathcal{O}(\mathrm{U}_{\sigma})=K\langle S_{\sigma}\rangle$,
$\mathcal{O}(\mathrm{U}_{\tau})=K\langle S_{\tau}\rangle$ and
$\mathcal{O}(\mathrm{U}_{\sigma\bigcap\tau})=K\langle
S_{\sigma\bigcap\tau}\rangle$, we only have to prove that
$K\langle S_{\sigma}\rangle$ and $K\langle S_{\tau}\rangle$ can
generate $K\langle S_{\sigma\bigcap\tau}\rangle$. It is obvious
that $S_{\sigma\bigcap\tau}$ is the sum of $S_{\sigma}$ and
$S_{\tau}$ due to the relation (\ref{sep}).
\end{proof}

The rigid spaces we get are also integral rigid spaces.

{\bf Theorem 2.} $\mathrm{R}(\Delta)$ is integral for a given fan
$\Delta$ in $\mathrm{N}_{\mathbb{R}}$.
\begin{proof}Since intersect of $\mathrm{U}_{\sigma}\:(\sigma\in\Delta)$ is not empty, we only have to prove that
$K\langle S_{\sigma}\rangle$ is integral, which is the next
proposition.
\end{proof}

{\bf Proposition 2.} If $\sigma$ is a strongly convex rational
polyhedral cone in $\mathrm{N}_{\mathbb{R}}$, $K\langle
S_{\sigma}\rangle$ is integral.
\begin{proof}
Choose a basis $u_{1},...,u_{n}$ of $\mathrm{M}_{\mathbb{R}}$,
which lie in $ \sigma^{\vee}\cap \mathrm{M}$. We can define an
order of $\sigma^{\vee}\cap \mathrm{M}$ compatible with the
additive operation in $\sigma^{\vee}\cap \mathrm{M}$ in the
following way. Let $u,v\in \sigma^{\vee}\cap \mathrm{M}$, and
write $u=\sum_{i}a_{i}u_{i},\:v=\sum_{i}b_{i}u_{i}$, with $\:
a_{i},\; b_{i}\in\mathbb{Q}$. We say $u<v$, if and only if there
exists an $i$,  $1\leq i\leq n$, such that $a_{j}=b_{j}$ for all
$1\leq j<i$ and $a_{i}<b_{i}$.

For any nonzero element $f$ in $K\langle S_{\sigma}\rangle$, write
$f=\sum_{u}c_{f,u}\chi^{u}$ with $c_{f,u}\in K$. Let $u_{f}\in
\sigma^{\vee}\cap \mathrm{M}$ be $\min\{u\in \sigma^{\vee}\cap
\mathrm{M};\; |c_{f,u}|\;\mathrm{is\; maximal}\}$.

For two nonzero elements $f,\;g$ in $K\langle S_{\sigma}\rangle$,
it is easy to see that $|c_{fg,u}|<|c_{f,u_{f}}||c_{g,u_{g}}|$ for
$u<u_{f}+u_{g}$, that $|c_{fg,u}|\leq |c_{f,u_{f}}||c_{g,u_{g}}|$
for $u>u_{f}+u_{g}$, and that
$|c_{fg,u_{f}+u_{g}}|=|c_{f,u_{f}}||c_{g,u_{g}}|$. Therefore $fg$
is not zero and $u_{fg}=u_{f}+u_{g}$. We have finished the proof.
\end{proof}

We write $pt$ for the fan defined by the origin. Then
\begin{equation}K\langle S_{pt}\rangle=K\langle X_{1},...,
X_{n},X_{1}^{-1},...,X_{n}^{-1}\rangle. \end{equation}
$\mathrm{U}_{pt}$ is a group rigid space (a torus) and has an
action on $\mathrm{R}(\Delta)$. In fact $K\langle S_{pt}\rangle$
is a Hopf affinoid algebra over $K$. Its comultiplication $m^{*}$
is defined by $m^{*}(X_{i})=X_{i}\otimes X_{i}$, its coinverse $s$
is defined by $s(X_{i})=X_{i}^{-1}$, and its counite $\varepsilon$
is defined by $\varepsilon(X_{i})=1$ for $1 \leq i \leq n$. For
any cone $\sigma$, there is a natural map $K\langle
S_{\sigma}\rangle\hookrightarrow K\langle S_{pt}\rangle$. It's
easy to see that $m^{*}(K\langle S_{\sigma}\rangle)\subseteq
K\langle S_{\sigma}\rangle\otimes K\langle S_{\sigma}\rangle
\subseteq K\langle S_{\sigma}\rangle \otimes K\langle
S_{pt}\rangle$, which defines a natural action of
$\mathrm{U}_{pt}$ on $\mathrm{U}_{\sigma}$. If $\tau$ is a face of
$\sigma$, the action of $\mathrm{U}_{pt}$ on $\mathrm{U}_{\tau}$
is the same as the restriction of the action of $\mathrm{U}_{pt}$
on $\mathrm{U}_{\sigma}$. Therefore, patching all these actions,
we get an action of $\mathrm{U}_{pt}$ on $\mathrm{R}(\Delta)$.

We point out that, since there is a natural affine covering over
$\mathrm{R}(\Delta)$, we can use it to calculate the Cech
cohomology for any given coherent sheaf over $\mathrm{R}(\Delta)$.

We can show that when the support of $\Delta$ is the whole
$\mathrm{N}_{\mathbb{R}}$, $\mathrm{R}(\Delta)$ is a proper rigid
space (see \cite{Fulton} for its proof).

For a given fan $\Delta$, we can also construct a toric scheme
$\mathrm{T}(\Delta)$ in the original sense, then its rigid
analytification (definition of rigid analytification, see
\cite{Bosch}) $\mathrm{RT}(\Delta)$ is also a rigid space, but it
is not isomorphic to $\mathrm{R}(\Delta)$ in general. For example,
for a simple cone $\sigma$, $\mathrm{RT}(\sigma)$ hasn't a finite
affinoid covering, however, $\mathrm{R}(\sigma)$ is an affinoid
space. But the authors don't know whether $\mathrm{R}(\Delta)$ and
$\mathrm{RT}(\Delta)$ are isomorphic when the support of $\Delta$
is the whole $\mathrm{N}_{\mathbb{R}}$.

\subsection{Examples}
In this subsection, we give two examples of toric rigid spaces.

 {\bf Example 1.} The projective rigid spaces $\mathbb{P}_{n}$ over $\bar{\mathbb{Q}}_{p}$.

 1. $n=2$

 Assume $\mathrm{N}=\mathbb{Z}\times\mathbb{Z}$ with a basis $\{u_{1},u_{2}\}$. We define three strongly convex rational
polyhedral cones $\sigma_{1}$,
 $\sigma_{2}$, $\sigma_{3}$ in
$\mathrm{N}_{\mathbb{R}}$ in the following way. $\sigma_{1}$ is
the cone generated by $u_{1}$ and $u_{2}$, $\sigma_{2}$ is the
cone generated by $-u_{1}-u_{2}$ and $u_{2}$, and $\sigma_{3}$ is
the cone generated by $-u_{1}-u_{2}$ and $u_{1}$. Let $\Delta$ be
the fan
$$\{\sigma_{1},\sigma_{2},\sigma_{3},\sigma_{1}\cap\sigma_{2},\sigma_{1}\cap\sigma_{3},\sigma_{2}\cap\sigma_{3},\sigma_{1}\cap\sigma_{2}\cap\sigma_{3}\}.$$
Write $X_{i}$ for $\chi^{u_{i}}$ ($i=1,2$) for simple. Then
$\bar{Q}_{p}\langle S_{\sigma_{1}}\rangle=\bar{Q}_{p}\langle
X_{1},X_{2}\rangle$, $\bar{Q}_{p}\langle S_{\sigma_{2}}\rangle
=\bar{Q}_{p}\langle X_{1}^{-1},X_{1}^{-1}X_{2}\rangle$, and
$\bar{Q}_{p}\langle S_{\sigma_{3}}\rangle =\bar{Q}_{p}\langle
X_{2}^{-1},X_{2}^{-1}X_{1}\rangle$. $\mathrm{R}(\Delta)$ is the
two dimensional projective rigid space $\mathbb{P}^{2}$ over
$\bar{\mathbb{Q}}_{p}$. $\mathrm{U}_{\sigma_{1}}$ corresponds to
$\{[x_{0},x_{1},x_{2}]\in \mathbb{P}_{2}:
|x_{1}|\leq|x_{0}|,|x_{2}|\leq|x_{0}|\}$,
$\mathrm{U}_{\sigma_{2}}$ corresponds to $\{[x_{0},x_{1},x_{2}]\in
\mathbb{P}_{2}: |x_{0}|\leq|x_{1}|,|x_{2}|\leq|x_{1}|\}$, and
$\mathrm{U}_{\sigma_{2}}$ corresponds to $\{[x_{0},x_{1},x_{2}]\in
\mathbb{P}_{2}: |x_{0}|\leq|x_{2}|,|x_{1}|\leq|x_{2}|\}$.

2. $n\geq 3$

It is similar to the case $n=2$. Assume
$\mathrm{N}=\mathbb{Z}^{\oplus n}$ with a basis
$\{u_{1},u_{2},\dots, u_{n}\}$. There are $n+1$ vectors
$v_{i}(0\leq i\leq n)$ with $v_{i}=u_{i}$ for $1\leq i\leq n$ and
$v_{0}=-u_{1}-u_{2}\dots-u_{n}$. For any proper subset $I$ of
$\{0,1,\dots ,n\}$, define
\begin{equation}
\sigma_{I}:\;=
\begin{cases}
\text{the cone of origin} & \text{if $I$ is empty},\\
\text{the cone generated by $v_{i}(i\in I)$}& \text{if $I$ is not
empty}.
\end{cases}
\end{equation}
And set
$$\triangle:=\{\sigma_{I}|\:I\subset \{0,1,\dots ,n\}\}.$$
Then $\mathrm{R}(\Delta)$ is the $n$ dimensional projective rigid
space $\mathbb{P}^{n}$.

$\mathrm{R}(\Delta)$ and $\mathrm{R}\mathrm{T}(\Delta)$ are the
same this time.

{\bf Example 2.} Hirzebrunch surfaces.

Let $a$ be a fixed positive integer. Assume
$\mathrm{N}=\mathbb{Z}\times\mathbb{Z}$ with a basis
$\{u_{1},u_{2}\}$. We define four strongly convex rational
polyhedral cones $\sigma_{1}$,
 $\sigma_{2}$, $\sigma_{3}$, $\sigma_{4}$ in
$\mathrm{N}_{\mathbb{R}}$ in the following way. $\sigma_{1}$ is
the cone generated by $u_{1}$ and $u_{2}$, $\sigma_{2}$ is the
cone generated by $u_{1}$ and $-u_{2}$, $\sigma_{3}$ is the cone
generated by $-u_{1}+a u_{2}$ and $-u_{2}$, and $\sigma_{4}$ is
the cone generated by $-u_{1}+a u_{2}$ and $u_{2}$. The four
corresponding affinoid spaces are $\bar{Q}_{p}\langle
S_{\sigma_{1}}\rangle =\bar{Q}_{p}\langle X_{1},X_{2}\rangle $,
$\bar{Q}_{p}\langle S_{\sigma_{2}}\rangle =\bar{Q}_{p}\langle
X_{1},X_{2}^{-1}\rangle$, $\bar{Q}_{p}\langle
S_{\sigma_{3}}\rangle =\bar{Q}_{p}\langle
X_{1}^{-1},X_{1}^{-a}X_{2}^{-1}\rangle $, $\bar{Q}_{p}\langle
S_{\sigma_{4}}\rangle =\bar{Q}_{p}\langle
X_{1}^{-1},X_{1}^{a}X_{2}\rangle $. Let $\Delta$ be the fan
{\small
$$\{\sigma_{1},\sigma_{2},\sigma_{3},\sigma_{4},\sigma_{1}\cap\sigma_{2},\sigma_{1}\cap\sigma_{4},
\sigma_{2}\cap\sigma_{3},\sigma_{3}\cap\sigma_{4},\sigma_{1}\cap\sigma_{2}\cap\sigma_{3}\cap\sigma_{4}
\}.$$} We get a rigid space $\mathrm{R}(\Delta)$ and call it a
Hirzebrunch surface.

Let $\mathrm{N}_{1}=\mathbb{Z}$ with a basis $\{v\}$. Let
$\tau_{1}=\{rv:\;r\geq 0\}$, $\tau_{2}=\{rv:\;r\leq 0\}$, and
$\Delta_{1}$ be the cone $\{\tau_{1},\tau_{2},\{0\}\}$. Then
$\mathrm{R}(\Delta_{1})$ is just the one dimensional projective
space. The  linear map
$\mathrm{N}_{\mathbb{R}}\rightarrow\mathrm{N}_{1\;\mathbb{R}}$
defined by $ru_{1}+su_{2}\mapsto rv$ induces a map from the fan
$\Delta_{1}$ to the fan $\Delta$, which determines a morphism
$\mathrm{R}(\Delta)\longrightarrow\mathrm{R}(\Delta_{1})$. This
morphism makes $\mathrm{R}(\Delta)$ a $\mathbb{P}^{1}$-bundle over
$\mathbb{P}^{1}$. Especially, when $a=0$, $\mathrm{R}(\Delta)$ is
the trivial bundle $\mathbb{P}^{1}\times\mathbb{P}^{1}$ over
$\mathbb{P}^{1}$.

\section{Reductions of toric rigid spaces}
Reduction theory is very important in rigid geometry. In this
section, we will study the reductions of toric rigid spaces. At
first, let's recall the reductions of affinoid algebras. For more
details of reductions, see \cite{BGR}.

Assume $K$ is a field with a non-Archimedean valuation. Let $A$ be
an affinoid algebra over $K$. There is a supremum (semi-)norm over
$A$ defined in the following way. For any $f\in A$, define
\begin{equation}
|f|_{sup}=\sup_{x\in \mathrm{Max} A}|f(x)|.
\end{equation}
Then define the reduction $\widetilde{A}$ of $A$ to be $\{ f\in A:
 |f|_{sup}\leq 1\}/\{f\in A:|f|_{sup}<1\}$. Reduction of
 $\mathrm{Sp}(A)$ is just defined to be
 $\mathrm{Spec}(\widetilde{A})$.

 We have the following lemma.

 {\bf Lemma 2.} Let $\sigma$ is a cone in $\mathrm{N}_{\mathbb{R}}$ and $S_{\sigma}$ be as in section 2. For $f\in K\langle
 S\rangle$, write
 $f=\sum\limits_{u\in\;\sigma}a_{u}\chi^{u}$
 with $a_{u}\in K$, $a_{u}\rightarrow 0$. Then
 \begin{equation}
 |f|_{sup}=\max_{u\in\;\sigma}|a_{u}|.
 \end{equation}
\begin{proof}
It's well known that, for any given valuation $|\cdot \:|$ over an
integral affinoid algebra $A$,
$|f|_{sup}=\lim\limits_{n\rightarrow
\infty}|f^{n}|^{\frac{1}{n}},\;\;\forall\;f\in A$. We use this
result to prove our assertion.

We write $|f|=\max\limits_{u\in\;\sigma}|a_{u}|$ for
$f=\sum\limits_{u\in\;\sigma}a_{u}\chi^{u}$ in $K\langle
 S\rangle$. It is easy to see that it defines a norm over $K\langle
 S\rangle$. From the proof of proposition 2, we can
 show that $|f|^{n}=|f^{n}|$. Therefore, we get
 $|f|_{sup}=\max_{u\in\;\sigma}|a_{u}|.$
\end{proof}

 Let $k$ be the residue field of $K$, then we have the following
 corollary.

 {\bf Corollary.} The reduction of $K\langle S_{\sigma}\rangle$ is
 $k[S_{\sigma}]$.

 For a separated rigid space $X$, its reduction always depends on a
 choice of an affinoid covering $\{U_{i}\}$. Open subsets $U_{i}$ and
 their intersections are all affinoid, so they have reductions
 defined as above. Patching them together, we get a reduction of
 $X$. It is a scheme over $k$. We call it the reduction of $X$ according to (the
 affinoid covering) $\{U_{i}\}$. In general, different affinoid
 coverings give different reductions.

 For a fan $\Delta$ in $\mathrm{N}$, we have defined a toric rigid space $\mathrm{R}(\Delta)$ in section 2.
There is a natural affinoid covering $\mathrm{Sp}(K\langle
S_{\sigma})\rangle\;(\sigma\in \Delta)$ of this rigid space
$\mathrm{R}(\Delta)$. Each $\mathrm{Sp}(K\langle
S_{\sigma}\rangle)$ has a reduction
$\mathrm{Spec}(k[S_{\sigma}])$. If $\sigma$ is a face of $\tau$,
the open immersion $\mathrm{Sp}(k\langle
S_{\sigma}\rangle)\hookrightarrow\mathrm{Sp}(k\langle
S_{\tau}\rangle)$ induces an open immersion
$\mathrm{Spec}(k[S_{\sigma}])\hookrightarrow\mathrm{Spec}(k[S_{\tau}])$.
Patching all the $\mathrm{Spec}(S_{\sigma})$ according to these
open immersions, we get a reduction of $\mathrm{R}(\Delta)$. It's
easy to see that this reduction is exactly the toric scheme
$\mathrm{T}(\Delta)_{k}$ over $k$.

For example, we calculate the reduction of the projective rigid
space $\mathbb{P}_{2}$ over $\bar{\mathbb{Q}}_{p}$ (example 1 in
section 2). It's patched by three spaces: {\scriptsize
$\mathrm{Spec}(\bar{\mathbb{F}}_{p}[X_{1},X_{2}])$,
$\mathrm{Spec}(\bar{\mathbb{F}}_{p}[X_{1}^{-1},X_{1}^{-1}X_{2}])$}
and
{\scriptsize$\mathrm{Spec}(\bar{\mathbb{F}}_{p}[X_{2}^{-1},X_{2}^{-1}X_{1}])$.}
Let $X_{1}=T_{1}/T_{0}, X_{2}=T_{2}/T_{0} $, these spaces are
corresponding to $\{[t_{0}:t_{1}:t_{2}]\in
\mathbb{P}_{2}|t_{0}\neq 0\}$, $\{[t_{0}:t_{1}:t_{2}]\in
\mathbb{P}_{2}|t_{1}\neq 0\}$, $\{[t_{0}:t_{1}:t_{2}]\in
\mathbb{P}_{2}|t_{2}\neq 0\}$. They form an affine covering of the
two dimensional projective spaces $\mathbb{P}_{2}$ over
$\bar{\mathbb{F}}_{p}$.

\end{document}